\documentclass[12pt, a4paper]{article}

\usepackage{tikz}
\usetikzlibrary{
    positioning,   
    arrows.meta,   
    shapes.geometric,
    calc
}
\usepackage{caption}
\usepackage{algorithm}
\usepackage{algpseudocode}

\usepackage[margin=0.9in]{geometry}
\usepackage[utf8]{inputenc} 
\usepackage{amsmath,amssymb,amsfonts}
\usepackage{graphicx}
\usepackage{xcolor}
\usepackage{natbib}
\usepackage{lscape}

\usepackage{titlesec}
\usepackage[normalem]{ulem}

\titlespacing{\section}{0pt}{6pt}{3pt}

\titlespacing{\subsection}{0pt}{3pt}{2pt} 
\titlespacing{\subsubsection}{0pt}{3pt}{2pt}

\usepackage{hyperref}
\usepackage{array}
\usepackage{lineno}
\usepackage{subfig}

\usepackage{float}
 \usepackage[toc,page]{appendix}
\usepackage{longtable}
\usepackage{hhline}
\usepackage{tabularray} 
\UseTblrLibrary{amsmath} 
\usepackage{amsthm}
\theoremstyle{definition}
\newtheorem{definition}{Definition}[section]
\newcounter{example}[section]

\theoremstyle{remark}
\newtheorem*{remark}{Remark}

\usepackage{arydshln}
\usepackage{diagbox}
\usepackage{rotating}
\usepackage{colortbl}
\hypersetup{
    colorlinks=true,
    linkcolor=blue,
    filecolor=magenta,      
    urlcolor=cyan,
}
\date{}
\usepackage{lipsum}

\usepackage{rotating}
\usepackage{cleveref}
\usepackage{multirow}
    \usepackage{hhline,booktabs}
    \usepackage{makecell}
\usepackage{epstopdf}
\begin{document}

\title{Neural-Guided Domain Restriction to Accelerate Pseudospectra
Computation for Structured Non-normal Banded Matrices}

\author{
Amit Punia$^{1}$
\and
Rakesh Kumar$^{2}$\thanks{Corresponding author. Email: \href{mailto:r.k.dhiman.math@gmail.com}{r.k.dhiman.math@gmail.com}}
\and
Madan Lal$^{3}$
}

\maketitle 

\begin{center}
{\small
$^{1}$Jai Narain Vyas University, Jodhpur, Rajasthan, India\\
Email: \href{mailto:amitpunia.maths@gmail.com}{amitpunia.maths@gmail.com}

\vspace{0.3em}

$^{2}$SRM University-AP, Amaravati, Andhra Pradesh, India\\
Email: \href{mailto:r.k.dhiman.math@gmail.com}{r.k.dhiman.math@gmail.com}

\vspace{0.3em}

$^{3}$Jai Narain Vyas University, Jodhpur, Rajasthan, India\\
Email: \href{mailto:madan.lakhani288@gmail.com}{madan.lakhani288@gmail.com}
}
\end{center}

\begin{abstract}
Computing pseudospectra of non-normal matrices is essential for understanding the stability and transient behavior of dynamical systems. Such analysis is critical in applications including fluid dynamics, control systems, and differential operators, where non-normality can lead to significant transient amplification and sensitivity to perturbations that are not captured by eigenvalue analysis alone.  At large scales, commonly used numerical approaches for pseudospectra computation can become computationally demanding, as they require repeated auxiliary computations to identify spectrally sensitive regions in the complex plane.

We present a neural network-based approach that predicts sensitive regions directly from matrix features, thereby avoiding exhaustive pseudospectra evaluation across the entire complex plane. We calibrate the prediction threshold on validation data to ensure reliable coverage of sensitive regions. The trained neural network guides the selection of grid points requiring full computation, enabling focused computation only where necessary. The approach provides a practical preprocessing strategy for efficient pseudospectra computation. Numerical experiments on non-normal banded matrices demonstrate substantial speedup compared to full grid-based numerical evaluation while maintaining high accuracy in identifying sensitive regions. 

\textbf{Keywords.} Pseudospectra, Non-normal matrices, Neural networks, 
Spectral analysis, Adaptive methods, Banded matrices

\textbf{AMS subject classification.} 15A18, 65F15, 68T07, 47A10

\end{abstract}

\newpage

\section{Introduction}

For a comprehensive understanding of physical models in science and engineering, the formulation of an appropriate mathematical model is essential. Valuable insights into system behavior can often be obtained from the eigenvalues or spectrum of the governing operator, highlighting their fundamental role in both theoretical and applied mathematics. In general, eigenvalues provide insight into algorithmic behavior, physical phenomena, and the qualitative structure of dynamical systems \cite{trefethen1997pseudospectra}. They are widely used across diverse fields including quantum mechanics, structural mechanics, fluid dynamics, numerical analysis, and probability theory, where they help characterize stability, asymptotic behavior, and rates of convergence or divergence.

Despite their importance, eigenvalues alone are often insufficient to describe the behavior of non-normal operators. Several models in fluid mechanics, such as the instability of magnetic plasmas, the formation of cyclones, and fluid flow in circular pipes, demonstrate the limitations of eigenvalue-based analysis \cite{trefethen2005spectra}. A classical example is the pipe Poiseuille flow problem, where the flow remains smooth at low velocities but transitions to turbulence at higher velocities \cite{trefethen1993hydrodynamic, schmid2001stability, schmid2007nonmodal}. In such cases, eigenvalue analysis fails to explain experimentally observed instabilities. In contrast, pseudospectral analysis captures transient growth and sensitivity effects, thereby providing a more accurate description of system behavior.

Pseudospectra also play a crucial role in numerical analysis. They provide valuable insight into stiffness, numerical instability of discretized differential equations, and the convergence behavior of iterative methods for non-symmetric matrix problems \cite{trefethen2005spectra, embree1999pseudospectra, trefethen1996finite, hairer1996solving}. These phenomena are closely linked to the non-normality of the underlying operators, where even small perturbations can lead to significant changes in system dynamics.

In general, non-normal and non-Hermitian problems exhibit behavior that cannot be fully characterized by eigenvalues alone. Such systems require a more refined analytical framework, where pseudospectra provide deeper insight into stability, sensitivity, and transient dynamics. Applications such as optical wave propagation further illustrate this sensitivity in complex media \cite{rotter2009light}.

Recent advances in machine learning (ML) and artificial intelligence (AI) have demonstrated strong performance in high-dimensional approximation, surrogate modeling, and data-driven scientific computing \cite{goodfellow2016deep, karniadakis2021physics, lu2021deeponet}. These approaches are particularly effective in problems where analytical characterization is challenging, such as non-self-adjoint eigenvalue problems in numerical linear algebra. Such problems are highly sensitive to perturbations, often leading to spectral pollution, transient growth, and numerical instabilities.

Unlike self-adjoint systems, which are supported by well-established theoretical frameworks such as Sturm--Liouville theory \cite{zettl2005sturm}, non-self-adjoint operators do not admit a similarly complete theory \cite{trefethen2005spectra}. This gap motivates the development of hybrid computational approaches that combine data-driven learning with classical numerical methods.

To address these challenges, we propose a scalable and interpretable hybrid framework for predicting pseudospectral behavior. The approach integrates neural network-based prediction with selective singular value computations to identify regions of spectral sensitivity efficiently. By leveraging the pattern recognition capabilities of neural networks, the proposed method reduces the computational burden associated with full-grid pseudospectra evaluation while preserving essential spectral features.

To generate accurate training data, we employ the exclusion-region singular-value method, which is closely related to established pseudospectra computation techniques based on selective singular value evaluation \cite{wright2002large, wright2002eigtool}. This approach ensures that the minimum singular value is computed exactly at relevant grid points, thereby providing reliable and physically consistent labels for training. Moreover, it preserves fine pseudospectral structures, such as sharp sensitivity regions and localized features near eigenvalues, which are essential for capturing subtle instability mechanisms in non-normal systems.

The computational efficiency of this approach enables the construction of large and diverse datasets consisting of thousands of random matrices and their corresponding pseudospectra. Such datasets are crucial for training robust machine learning models and mitigating overfitting. To the best of our knowledge, relatively limited work exists on directly predicting matrix pseudospectra using deep learning. However, recent advances in scientific machine learning and operator learning suggest that data-driven approaches can effectively approximate complex operator-dependent mappings \cite{karniadakis2021physics, lu2021deeponet}.

\textbf{Contributions of this work are as follows:}
\begin{itemize}
\item We propose a hybrid data-driven framework for efficient approximation of pseudospectra in non-normal matrices.
\item We develop a neural network-based predictor to identify regions of spectral sensitivity, significantly reducing the number of required singular value computations.
\item We design a data generation strategy based on selective singular value evaluation to produce accurate and physically consistent training labels.
\item We introduce a hierarchical coarse-to-fine prediction strategy to reduce the number of neural network evaluations.
\item We demonstrate substantial computational speedup while maintaining high accuracy through extensive numerical experiments.
\end{itemize}

The remainder of this paper is organized as follows. Section~\ref{sec:preliminaries} introduces the mathematical foundations of pseudospectra, Section~\ref{sec:neural-approach} presents the proposed neural network framework, Section~\ref{sec:hybrid_approach} details the hybrid computational strategy, 
and Section~\ref{sec:experiments} reports numerical experiments and performance evaluation. Finally, Sections~\ref{sec:discussion} and~\ref{sec:conclusion} provide discussion and concluding remarks.

\section{Preliminaries}\label{sec:preliminaries}

For a matrix $\mathbf{A} \in \mathbb{C}^{n \times n}$, the \textit{pseudospectra} characterize the behavior of linear systems beyond what eigenvalues alone reveal. While eigenvalues determine asymptotic stability, pseudospectra govern transient dynamics and sensitivity to perturbations \cite{trefethen2005spectra}. This is particularly critical for non-normal matrices, where $\mathbf{A}\mathbf{A}^* \neq \mathbf{A}^*\mathbf{A}$, which arise in fluid dynamics \cite{schmid2001stability}, control theory \cite{trefethen1993hydrodynamic, hinrichsen2005mathematical}, and weather prediction \cite{farrell1996generalized}.

For non-normal matrices, the eigenvectors may fail to form a complete or orthogonal basis. Moreover, the eigenvector matrix can be highly ill-conditioned, indicating strong sensitivity of the system to even very small perturbations. 

\begin{figure}[h]
 \centering
 \includegraphics[width=.6\linewidth]{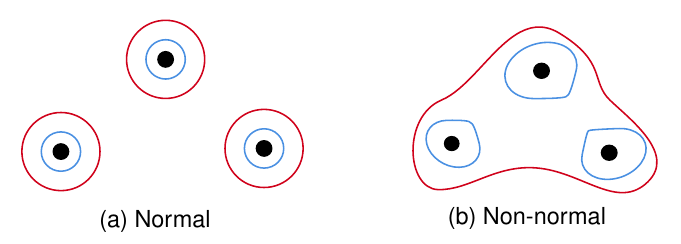}
 \caption{Effect of perturbation on normal and non-normal matrix eigenvalues}
 \label{fig:test1}
\end{figure}

As illustrated in Fig.~\ref{fig:test1}, for normal operators, the $\epsilon$-pseudospectrum $\Lambda_\epsilon(A)$ consists of points in $\mathbb{C}$ that lie within a distance $\epsilon$ of the spectrum $\sigma(A)$. In contrast, for non-normal operators, $\Lambda_\epsilon(A)$ can extend far beyond the eigenvalue locations, reflecting increased spectral sensitivity.

\begin{figure}[h]
 \centering
 \includegraphics[width=0.42\linewidth]{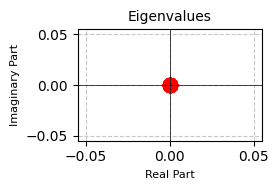}
 \hfill
 \includegraphics[width=0.42\linewidth]{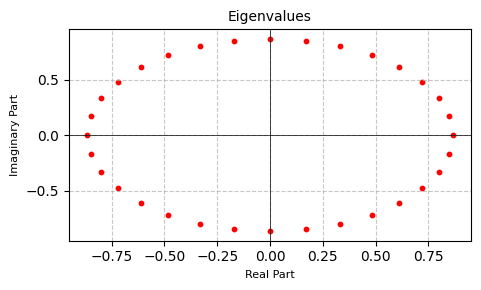}
 \caption{Effect of a perturbation on the eigenvalues of matrix $A$}
 \label{fig:32}
\end{figure}

\begin{figure}[h]
 \centering
 \includegraphics[width=0.5\linewidth]{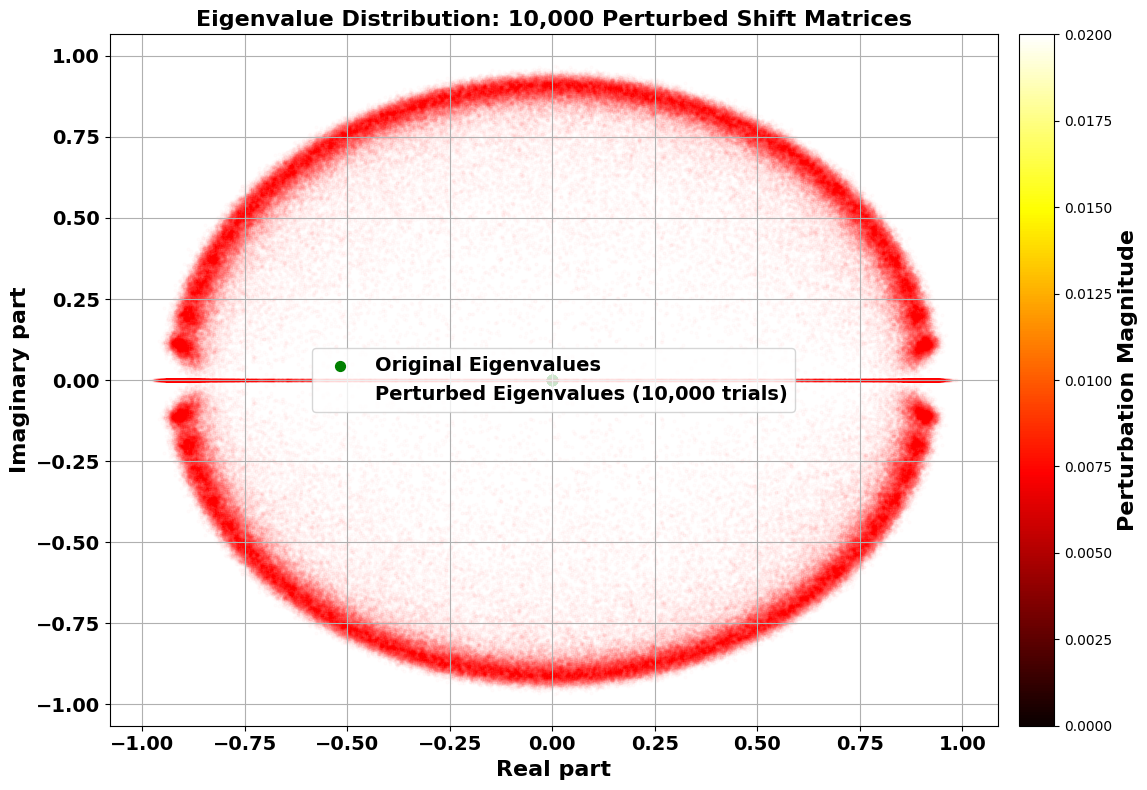}
 \caption{Distribution of eigenvalues under random perturbations}
 \label{fig:my_label1}
\end{figure}

Consider the Jordan block matrix $A$ of order $32$ with ones on the subdiagonal. All its eigenvalues are zero. However, even a small perturbation in a single entry of the first row of $A$ drastically affects the eigenvalues, as illustrated in Fig.~\ref{fig:32}. The eigenvalues of the perturbed matrix $B$ spread away from the origin and lie approximately on a circle centered at zero. Furthermore, Fig.~\ref{fig:my_label1} shows the distribution of eigenvalues under ten thousand random perturbations of the entries of $A$.

\[
A =
\begin{bmatrix}
0 & 0 & 0 & 0 & \textcolor{blue}{0} \\
1 & 0 & 0 & 0 & 0 \\
0 & 1 & 0 & 0 & 0 \\
0 & 0 & 1 & 0 & 0 \\
0 & 0 & 0 & 1 & 0
\end{bmatrix}_{32}
\hspace{2cm}
B =
\begin{bmatrix}
0 & 0 & 0 & 0 & \textcolor{red}{0.01} \\
1 & 0 & 0 & 0 & 0 \\
0 & 1 & 0 & 0 & 0 \\
0 & 0 & 1 & 0 & 0 \\
0 & 0 & 0 & 1 & 0
\end{bmatrix}_{32}
\]

These observations demonstrate that eigenvalues and eigenvectors alone are insufficient for analyzing non-normal matrices. The pseudospectrum, or $\epsilon$-spectrum, is defined as the set of all eigenvalues of $\mathbf{A}$ and those of nearby perturbed matrices \cite{trefethen2005spectra}. When dealing with highly non-normal problems, pseudospectral analysis is more reliable, as eigenvalue-based analysis may be insufficient.

Consider the dynamical system:
\begin{equation}
\frac{d\mathbf{x}}{dt} = \mathbf{A}\mathbf{x}.
\end{equation}

Even when all eigenvalues of $\mathbf{A}$ have negative real parts, non-normality can cause large transient growth before eventual decay. Classical eigenvalue analysis fails to capture this behavior, motivating the use of pseudospectral analysis.

\begin{definition}[Resolvent Norm]\label{def:rdist}
For a matrix $\mathbf{A}$, the resolvent is defined as $(z\mathbf{I} - \mathbf{A})^{-1}$ for $z \notin \sigma(\mathbf{A})$. If $\mathbf{A}$ is normal, then
\[
\|(z\mathbf{I} - \mathbf{A})^{-1}\| = \frac{1}{\mathrm{dist}(z, \sigma(\mathbf{A}))}.
\]
\end{definition}

\begin{remark}
For non-normal matrices, the resolvent norm may significantly deviate from $\frac{1}{\mathrm{dist}(z, \sigma(\mathbf{A}))}$, reflecting high spectral sensitivity.
\end{remark}

\begin{definition}[Pseudospectrum]\cite{trefethen1997pseudospectra,trefethen2005spectra}
For $\mathbf{A} \in \mathbb{C}^{n \times n}$ and $\epsilon > 0$, the $\epsilon$-pseudospectrum is defined as
\begin{equation}\label{eq:pseudospectrum}
\begin{aligned}
\Lambda_\epsilon(\mathbf{A}) 
&= \left\{ z \in \mathbb{C} : \|(z\mathbf{I} - \mathbf{A})^{-1}\| \geq \frac{1}{\epsilon} \right\} \\
&= \bigcup_{\|\mathbf{E}\| \leq \epsilon} \lambda(\mathbf{A} + \mathbf{E}) \\
&= \left\{ z \in \mathbb{C} : \sigma_{\min}(z\mathbf{I} - \mathbf{A}) \leq \epsilon \right\}.
\end{aligned}
\end{equation}
Here, $\lambda(\cdot)$ denotes the spectrum and $\sigma_{\min}(\cdot)$ denotes the smallest singular value.
\end{definition}

The second characterization shows that $\Lambda_\epsilon(\mathbf{A})$ contains all eigenvalues of matrices within an $\epsilon$-neighborhood of $\mathbf{A}$, thereby quantifying spectral sensitivity.

\begin{remark}[Non-normality and Spectral Sensitivity]
For normal matrices, the pseudospectra consist of $\epsilon$-disks centered at eigenvalues. For non-normal matrices, $\Lambda_\epsilon(\mathbf{A})$ can extend far beyond the eigenvalue locations, indicating strong sensitivity to perturbations \cite{trefethen1997pseudospectra, schmid2001stability}.
\end{remark}

\subsection{Singular Value Decomposition}

The standard approach to computing pseudospectra discretizes the complex plane on a rectangular grid 
\[
\mathcal{G} = \{z_{ij} : i = 1, \ldots, N_y,\; j = 1, \ldots, N_x\},
\]
where
\[
z_{ij} = x_j + \mathrm{i} y_i, \quad x_j \in [a, b], \quad y_i \in [c, d],
\]
for a domain $[a, b] \times [c, d] \subset \mathbb{C}$ covering the region of interest.

For each grid point $z \in \mathcal{G}$, we compute
\[
\sigma_{\min}(z) = \sigma_{\min}(zI - A)
\]
via the singular value decomposition (SVD), which forms the basis of standard pseudospectra computation methods 
\cite{trefethen2005spectra, wright2002large, braconnier1996lanczos, lui1997computation}.

If $zI - A = U \Sigma V^*$ is the SVD, where $\Sigma = \mathrm{diag}(\sigma_1, \ldots, \sigma_n)$ with $\sigma_1 \geq \cdots \geq \sigma_n \geq 0$, then $\sigma_{\min}(z) = \sigma_n$.

The region of interest, referred to as the \emph{sensitive zone}, is defined as
\begin{equation}
\label{eq:sensitive_zone}
\mathcal{S}_\varepsilon = \left\{ z \in \mathcal{G} : \sigma_{\min}(z) \leq \varepsilon \right\}.
\end{equation}

Visualization typically uses contour plots of $\log_{10}(\sigma_{\min}(z))$, where the level curve $\log_{10}(\varepsilon)$ approximates the pseudospectrum boundary.

\subsection{Computational Cost Analysis}

The primary computational bottleneck in pseudospectra computation arises from repeated singular value decomposition (SVD) evaluations, as observed in large-scale algorithms \cite{wright2002large}. For a matrix $A \in \mathbb{C}^{n \times n}$ and a grid $\mathcal{G}$ with $N = N_x \times N_y$ points, the total computational cost is
\begin{equation}
\label{eq:full_cost}
\mathcal{C}_{\text{full}} = N \cdot \mathcal{C}_{\text{SVD}}(n),
\end{equation}
where $\mathcal{C}_{\text{SVD}}(n)$ denotes the cost of a single SVD.

Using standard algorithms (e.g., Golub--Reinsch or divide-and-conquer methods), the cost scales as $\mathcal{C}_{\text{SVD}}(n) = O(n^3)$ floating-point operations \cite{golub2013matrix}.

For typical grid sizes ($N = 10^4$) and moderate matrix dimensions (e.g., $n = 64$), the total computational cost is on the order of $10^9$ floating-point operations, which becomes prohibitive for larger problems.

A key observation is that the sensitive region $\mathcal{S}_\varepsilon$ typically occupies only a small fraction of the grid. Let $\rho = |\mathcal{S}_\varepsilon| / N$. If this region could be identified \emph{a priori}, the computational cost could be reduced by a factor of
\begin{equation}
\label{eq:potential_speedup}
\text{Speedup}_{\text{potential}} = \rho^{-1}.
\end{equation}

In practice, exact identification of $\mathcal{S}_\varepsilon$ is not possible without computing the pseudospectra. The central challenge is therefore to design an efficient predictor that identifies $\mathcal{S}_\varepsilon$ with high recall while maintaining low computational overhead.

\section{Neural Network Approach}\label{sec:neural-approach}

In this section, the proposed neural network-based method for accelerating pseudospectra computation is discussed. A classifier is trained that predicts which regions of the complex plane require detailed singular value decomposition. This enables the restriction of expensive computations to a small fraction of the grid while maintaining high coverage of true sensitive zones. We now describe the matrix-level features used to characterize pseudospectral geometry and to condition the neural classifier.

\subsection{Matrix Feature Extraction}

To enable the neural network to learn geometric properties of pseudospectra from matrix structure alone, we extract a set of size-independent features motivated by classical matrix analysis, conditioning theory, and non-normality measures \cite{higham2002accuracy, trefethen2005spectra,trefethen1997pseudospectra} from input matrix $A \in \mathbb{R}^{n \times n}$.

To extract informative matrix features, we first normalize the input matrix as
\[
\tilde{A} = \frac{A}{\|A\|_F + \epsilon},
\]
where $\|A\|_F$ denotes the Frobenius norm and $\epsilon \approx 10^{-12}$ ensures numerical stability. This normalization enforces scale invariance, allowing the extracted features to remain comparable across matrices of different magnitudes.

Let $\{\lambda_i\}_{i=1}^n$ and $\{\sigma_i\}_{i=1}^n$ denote the eigenvalues and singular values of $A$, respectively, with singular values ordered as $\sigma_1 \geq \sigma_2 \geq \cdots \geq \sigma_n$. Based on these quantities, we compute a set of $23$ global features that capture key spectral, structural, and conditioning properties of the matrix, as summarized in Table~\ref{tab:matrix}. These features are designed to provide a compact representation of matrix characteristics that influence the geometry and extent of the pseudospectrum \cite{trefethen2005spectra, higham2002accuracy}.

\begin{table}[h!]
\begin{center}
{\renewcommand{\arraystretch}{1.7}
\caption{Global matrix features extracted from $A \in \mathbb{R}^{n \times n}$, including eigenvalue statistics, spectral spread, non-normality measures, conditioning, matrix norms, and sparsity descriptors. These features provide a compact representation of structural and spectral properties that influence the geometry of the pseudospectrum.} \label{tab:matrix}
\begin{tabular}{ |p{3.7cm}|p{3.6cm}|p{3.6cm}| p{3.6cm}|  }
\hline
\multicolumn{4}{|c|}{\textbf{Eigenvalue} } \\
\hline
$f_1 = \text{mean}(\mathfrak{Re}{(\lambda_i))}$ & $f_2 = \text{std}(\mathfrak{Re}{(\lambda_i)})$ & $f_3 = \min_i \mathfrak{Re}{(\lambda_i)}$& $f_4 = \max_i \mathfrak{Re}{(\lambda_i)}$ \\
$f_5 = \text{mean}(\mathfrak{Im}{(\lambda_i)})$ & $f_6 = \text{std}(\mathfrak{Im}{(\lambda_i)})$   & $f_7 = \min_i \mathfrak{Im}{(\lambda_i)}$ &$f_8 = \max_i \mathfrak{Im}{(\lambda_i)},$\\
\hline
\multicolumn{4}{|c|}{\textbf{Spectral spread} $\|$ \textbf{Non-normality measures}} \\
\hline
$f_9 = \max_i |\lambda_i|$ & $ f_{10} = \min_i |\lambda_i|$ & $f_{11} = \frac{ 1} { \|A\|_F}\|A - A^T\|_F$& $f_{12} = \frac{ 1}{ \|A\|_F}\|A - A^H\|_F$ \\
\hline
\multicolumn{4}{|c|}{\textbf{Conditioning} } \\
\hline
\multicolumn{4}{|c|}{$f_{13} = \log_{10}(\kappa(A) + \epsilon), \quad \text{Here, } \kappa(A) = \frac{\sigma_1}{\sigma_n},$ and $ \epsilon = 10^{-12} $ ensure numerical stability.} \\
\hline
\multicolumn{4}{|c|}{\textbf{Matrix norms}} \\
\hline
\multicolumn{4}{|c|}{$f_{14} = \frac{1 }{ \|A\|_F}\|A\|_2,$ \quad 
$f_{15} = \frac{1 }{ \|A\|_F}\|A\|_1,$ \quad 
$f_{16} = \frac{1}{ \|A\|_F}\|A\|_\infty .$} \\
\hline
\multicolumn{4}{|c|}{\textbf{Diagonal} $\&$ \textbf{off-diagonal properties}} \\
\hline
$f_{17} = \text{mean}(|\text{diag}(A)|)$ & $f_{18} = \text{std}(|\text{diag}(A)|)$ & $f_{19} = \text{mean}(|A - \text{diag}(\text{diag}(A))|)$ & $f_{20} = \text{std}(|A - \text{diag}(\text{diag}(A))|)$ \\
\hline
\multicolumn{4}{|c|}{\textbf{Sparsity $\&$ distribution}} \\
\hline
\multicolumn{4}{|c|}{$f_{21} = \frac{1}{n^2}\sum_{i,j} \mathbb{I}(|A_{ij}| > 10^{-10}), $\quad 
$f_{22} = \text{mean}(\tilde{A}^2),$ \quad 
$f_{23} = \text{std}(\tilde{A}^2).$} \\
\hline
\end{tabular}}
\end{center}
\end{table}

Given our focus on highly sensitive (non-normal) matrices, the associated eigenvectors are often ill-conditioned, leading to pseudospectra that extend significantly beyond the eigenvalue locations \cite{trefethen2005spectra}. Consequently, the global matrix features listed in Table~\ref{tab:matrix}, while informative, are not sufficient to fully capture the spatial extent and geometry of the pseudospectrum.

To address this limitation, we augment these features with an additional set of seven descriptors (see Tables~\ref{tab:matrixAdd} and \ref{tab:resolvent}) specifically designed to quantify spectral spread, non-normality, and resolvent growth. These additional features provide indirect but informative estimates of how far the pseudospectrum extends from the eigenvalues, thereby enabling the model to better predict sensitive regions in the complex plane.

\begin{table}[h!]
\begin{center}
{\renewcommand{\arraystretch}{1.7}
\caption{Additional matrix features designed to capture spectral sensitivity and non-normality effects. These include eigenvector conditioning, logarithmic non-normality ratios, and measures of singular value and eigenvalue spread, which are critical for estimating pseudospectral expansion beyond eigenvalue locations.} \label{tab:matrixAdd}
\begin{tabular}{ |p{3.7cm}|p{3.6cm}|p{3.6cm}| p{3.6cm}|  }
\hline
\multicolumn{4}{|c|}{\textbf{eigenvector condition
number $\&$ logarithmic non-normality ratio} } \\
\hline
\multicolumn{4}{|c|}{$f_{24} = \log_{10}(\kappa(V) + \delta), \quad f_{25} = \log_{10}\!\left( \frac{\|A - A^H\|_F}{\|A\|_F} + \epsilon \right).$ 

} \\
\hline
\multicolumn{4}{|c|}{\textbf{Singular value $\&$ Eigenvalue spread} } \\
\hline
\multicolumn{4}{|c|}{$f_{26} = \frac{\sigma_1 - \sigma_n}{\sigma_1 + \epsilon}, \quad f_{27} = \frac{\max_i |\lambda_i| - \min_i |\lambda_i|}{\max_i |\lambda_i| + \epsilon}$} \\
\hline
\end{tabular}}
\end{center}
\end{table}

Eigenvector conditioning plays a central role in spectral sensitivity, as characterized by classical perturbation results such as the Bauer-Fike theorem \cite{bauer1960norms}.

\begin{table}[h!]
\begin{center}
{\renewcommand{\arraystretch}{1.7}
\caption{Resolvent-based features computed at shifted points $z_\delta = \bar{\lambda} + \delta$ for $\delta \in \{0.5, 1.0, 2.0\}$, where $\bar{\lambda}$ denotes the eigenvalue centroid. These quantities approximate resolvent norm growth and provide insight into pseudospectral sensitivity away from the spectrum.}\label{tab:resolvent}
\begin{tabular}{ |p{5cm}|p{5cm}|p{5cm}|   }
\hline
\multicolumn{3}{|c|}{\textbf{Resolvent norm $\|(z_\delta I - A)^{-1}\| $ at distinct $\delta = 0.5,~1.0, ~2.0$} } \\
\hline
$f_{28} = \log_{10}(\|(z_{0.5} I - A)^{-1}\| ) $ & $f_{29} = \log_{10}(\|(z_{1.0} I - A)^{-1}\| )$ & $f_{30} = \log_{10}(\|(z_{2.0} I - A)^{-1}\| )$ \\
\hline
\end{tabular}}
\end{center}
\end{table}

Here, in Table~\ref{tab:resolvent}, $z_\delta = \bar{\lambda} + \delta$, where $\bar{\lambda}$ denotes the eigenvalue centroid. The resolvent norm $\|(z_\delta I - A)^{-1}\|$ is estimated by solving 
\[
(z_\delta I - A)x = b,
\]
for a random vector $b \sim \mathcal{N}(0, I)$ and computing the ratio $\|x\|_2 / \|b\|_2$.

These additional features complement global matrix descriptors by capturing spectral sensitivity and resolvent growth effects. In particular, non-normal matrices exhibit pseudospectral regions that can extend far beyond the eigenvalue locations, and this behavior depends on both eigenvector conditioning and resolvent growth \cite{trefethen2005spectra,trefethen1997pseudospectra}. Therefore, augmenting global descriptors with features related to spectral spread and resolvent estimates improves the predictive capability of the model.

The behavior of the pseudospectra is influenced by local geometry of the eigenvalues (Def.~\ref{def:rdist}). For each grid point $z = x + iy \in \mathbb{C}$, we compute three per-point eigenvalue distance features:
\begin{equation}\label{eq:dist}
g_1(z) = \min_{i} |z - \lambda_i|, \quad
g_2(z) = |z - \lambda^*|, \quad
g_3(z) = \frac{1}{n}\sum_{i=1}^n |z - \lambda_i|.
\end{equation}

These distance-based features provide local geometric context relative to the spectrum, enabling the network to adapt its predictions based on proximity to eigenvalues in the complex plane.
\begin{equation}\label{eq:f33}
\mathbf{f}(z, A) = [f_1,~f_2,~f_3,~ \ldots, f_{29},~ f_{30},~ g_1(z), ~g_2(z),~ g_3(z)] \in \mathbb{R}^{33}.
\end{equation}

These features are incorporated into the neural network architecture as described in the next subsection.

\subsection{Network Architecture}
In our approach, the classifier consists of two parallel pathways that process spatial coordinates and matrix features separately before fusion. We first describe the encoding used for spatial coordinates in the complex plane.

\subsubsection{Fourier Feature Encoding}

Grid coordinates $z = x + iy$ are first encoded using Fourier features to capture
periodic structure arising from oscillatory level sets of the resolvent norm in
the complex plane. The number of frequency bands is chosen to balance expressive
power with computational efficiency, providing sufficient resolution to capture
fine-scale pseudospectral structure without over-parameterization. For coordinate
input $\mathbf{c} = [x, y]^T \in \mathbb{R}^2$, we compute
\begin{equation}
\phi(\mathbf{c}) = \left[\mathbf{c}, \sin(2\mathbf{c}), \cos(2\mathbf{c}), 
\sin(4\mathbf{c}), \cos(4\mathbf{c}), \ldots, 
\sin(2^6 \mathbf{c}), \cos(2^6 \mathbf{c})\right]^T
\end{equation}
where the maximum frequency is $2^6 = 64$, the sine and cosine functions are applied element-wise across six frequency bands $(2^1,2^2, \ldots,2^6)$, yielding a compact multi-scale representation of the
spatial coordinates \cite{tancik2020fourier}.

The encoded spatial coordinates are then combined with matrix-dependent features
within a dual-path network architecture, which we describe next.

\subsubsection{Dual-Path Architecture} 

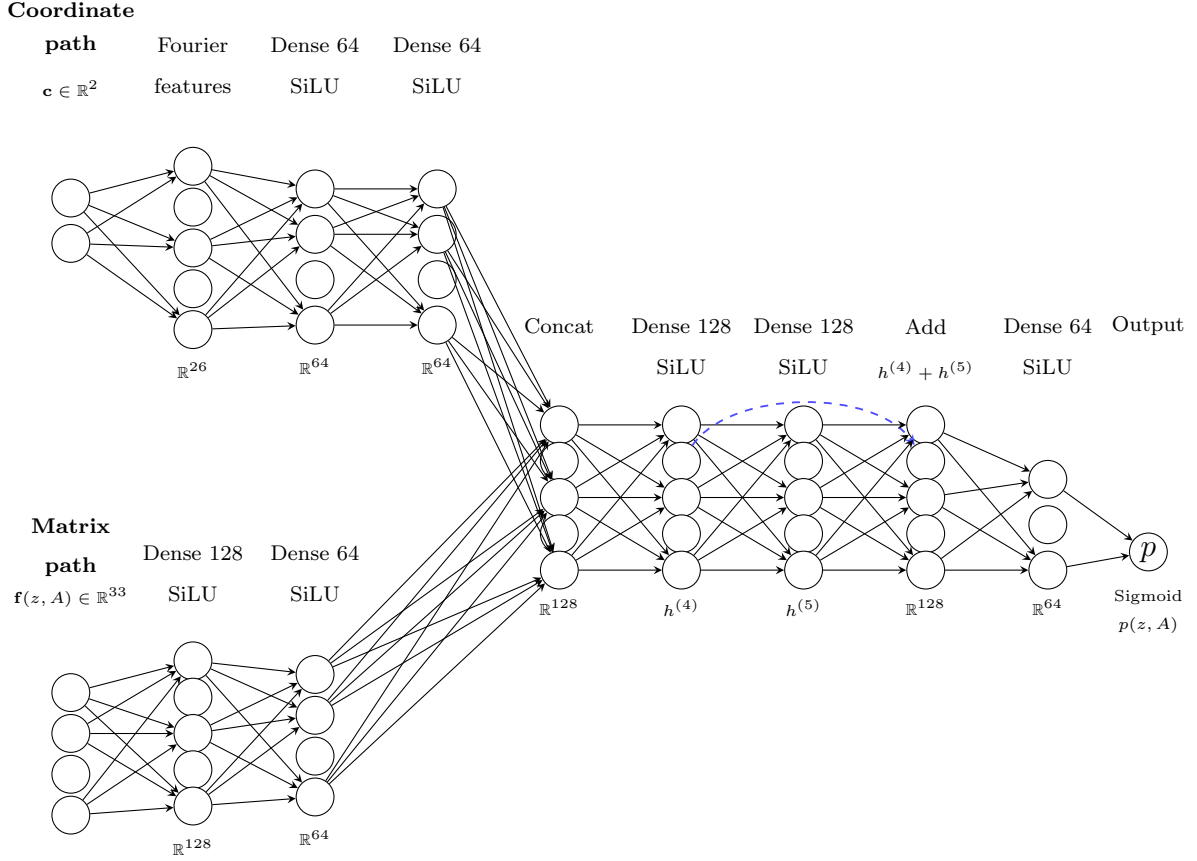
\begin{figure*}[ht]
\centering
\begin{tikzpicture}[
    xscale=0.95,
    yscale=1.2,
    neuron/.style={circle, draw, minimum size=5mm, inner sep=0pt, fill=white},
    label/.style={font=\scriptsize, align=center},
    dim/.style={font=\tiny},
    >=stealth
]

\node[label, yshift=6pt] at (0,4.2) {\textbf{Coordinate}};
\node[label, yshift=3pt] at (0,3.9) {\textbf{path}};
\node[dim] at (0,3.5) {$\mathbf{c}\in\mathbb{R}^2$};

\foreach \i in {1,2}
    \node[neuron] (z\i) at (0,2.8-0.5*\i) {};

\node[label, yshift=3pt] at (1.7,3.9) {Fourier};
\node[label, yshift=-2pt] at (1.7,3.6) {features};
\foreach \i in {1,...,5}
    \node[neuron] (ff\i) at (1.7,3.1-0.45*\i) {};
\node[dim] at (1.7,0.4) {$\mathbb{R}^{26}$};

\node[label, yshift=3pt] at (3.4,3.9) {Dense 64};
\node[label, yshift=-2pt] at (3.4,3.6) {SiLU};
\foreach \i in {1,...,4}
    \node[neuron] (c1\i) at (3.4,2.9-0.5*\i) {};
\node[dim] at (3.4,0.45) {$\mathbb{R}^{64}$};

\node[label, yshift=3pt] at (5.1,3.9) {Dense 64};
\node[label, yshift=-2pt] at (5.1,3.6) {SiLU};
\foreach \i in {1,...,4}
    \node[neuron] (c2\i) at (5.1,2.9-0.5*\i) {};
\node[dim] at (5.1,0.45) {$\mathbb{R}^{64}$};

\node[label, yshift=3pt] at (0,-1.4) {\textbf{Matrix}};
\node[label, yshift=-2pt] at (0,-1.7) {\textbf{path}};
\node[dim] at (0,-2.1) {$\mathbf{f}(z,A) \in\mathbb{R}^{33} $};

\foreach \i in {1,...,4}
    \node[neuron] (m\i) at (0,-2.7-0.45*\i) {};

\node[label, yshift=3pt] at (1.7,-1.7) {Dense 128};
\node[label, yshift=-2pt] at (1.7,-2.0) {SiLU};
\foreach \i in {1,...,5}
    \node[neuron] (m1\i) at (1.7,-2.4-0.4*\i) {};
\node[dim] at (1.7,-4.85) {$\mathbb{R}^{128}$};

\node[label, yshift=3pt] at (3.4,-1.7) {Dense 64};
\node[label, yshift=-2pt] at (3.4,-2.0) {SiLU};
\foreach \i in {1,...,4}
    \node[neuron] (m2\i) at (3.4,-2.5-0.45*\i) {};
\node[dim] at (3.4,-4.75) {$\mathbb{R}^{64}$};

\node[label, yshift=3pt] at (6.8,0.8) {Concat};
\foreach \i in {1,...,5}
    \node[neuron] (f\i) at (6.8,0.2-0.4*\i) {};
\node[dim] at (6.8,-2.21) {$\mathbb{R}^{128}$};

\node[label, yshift=3pt] at (8.5,0.8) {Dense 128};
\node[label, yshift=-2pt] at (8.5,0.5) {SiLU};
\foreach \i in {1,...,5}
    \node[neuron] (r1\i) at (8.5,0.2-0.4*\i) {};
\node[dim] at (8.5,-2.21) {$h^{(4)}$};

\node[label, yshift=3pt] at (10.2,0.8) {Dense 128};
\node[label, yshift=-2pt] at (10.2,0.5) {SiLU};
\foreach \i in {1,...,5}
    \node[neuron] (r2\i) at (10.2,0.2-0.4*\i) {};
\node[dim] at (10.2,-2.21) {$h^{(5)}$};

\node[label, yshift=3pt] at (11.9,0.8) {Add};
\node[dim] at (11.9,0.4) {$h^{(4)}+h^{(5)}$};
\foreach \i in {1,...,5}
    \node[neuron] (radd\i) at (11.9,0.2-0.4*\i) {};
\node[dim] at (11.9,-2.21) {$\mathbb{R}^{128}$};

\node[label, yshift=3pt] at (13.6,0.8) {Dense 64};
\node[label, yshift=-2pt] at (13.6,0.5) {SiLU};
\foreach \i in {1,...,3}
    \node[neuron] (r3\i) at (13.6,-0.3-0.5*\i) {};
\node[dim] at (13.6,-2.21) {$\mathbb{R}^{64}$};

\node[label, yshift=3pt] at (15.0,0.8) {Output};
\node[neuron] (out) at (15.0,-1.6) {$p$};
\node[dim] at (15.0,-2.1) {Sigmoid};
\node[dim] at (15.0,-2.4) {$p(z,A)$};

\foreach \i in {1,2}
    \foreach \j in {1,3,5}
        \draw[->] (z\i) -- (ff\j);

\foreach \i in {1,3,5}
    \foreach \j in {1,2,4}
        \draw[->] (ff\i) -- (c1\j);

\foreach \i in {1,2,4}
    \foreach \j in {1,2,4}
        \draw[->] (c1\i) -- (c2\j);

\foreach \i in {1,2,4}
    \foreach \j in {1,3,5}
        \draw[->] (m\i) -- (m1\j);

\foreach \i in {1,3,5}
    \foreach \j in {1,2,4}
        \draw[->] (m1\i) -- (m2\j);

\foreach \i in {1,2,4}
    \foreach \j in {1,3,5}
        \draw[->] (c2\i) -- (f\j);

\foreach \i in {1,2,4}
    \foreach \j in {1,3,5}
        \draw[->] (m2\i) -- (f\j);

\foreach \i in {1,3,5}
    \foreach \j in {1,3,5}
        \draw[->] (f\i) -- (r1\j);

\foreach \i in {1,3,5}
    \foreach \j in {1,3,5}
        \draw[->] (r1\i) -- (r2\j);

\draw[->, very thick, blue!70, dashed, line width=0.7pt] 
    (r12) .. controls +(0.7,0.8) and +(-0.7,0.8) .. (radd2);

\foreach \i in {1,3,5}
    \foreach \j in {1,3,5}
        \draw[->] (r2\i) -- (radd\j);

\foreach \i in {1,3,5}
    \foreach \j in {1,3}
        \draw[->] (radd\i) -- (r3\j);

\foreach \i in {1,3}
    \draw[->] (r3\i) -- (out);

\end{tikzpicture}
\caption{Neural network architecture for pseudospectra prediction. The coordinate pathway (top) processes spatial locations via Fourier features, while the matrix pathway (bottom) processes 33 combined features (30 global matrix properties + 3 per-point eigenvalue distances). Both pathways converge at the central concatenation layer. After concatenation, a residual block (with skip connection in blue) refines the representation before sigmoid output $p(z,A)=P(\text{sensitive}\mid z,A)$.}
\label{fig:architecture}
\end{figure*}

The Fourier-encoded coordinates are mapped through two fully-connected layers
\begin{align}
\mathbf{h}_c^{(1)} &= \text{SiLU}\left(W_c^{(1)} \phi(\mathbf{c}) + b_c^{(1)}\right), 
\quad W_c^{(1)} \in \mathbb{R}^{64 \times 26}, \\
\mathbf{h}_c^{(2)} &= \text{SiLU}\left(W_c^{(2)} \mathbf{h}_c^{(1)} + b_c^{(2)}\right), 
\quad W_c^{(2)} \in \mathbb{R}^{64 \times 64}.
\end{align}

The 33-dimensional feature vector $\mathbf{f}(z,A)$ is mapped through
\begin{align}
\mathbf{h}_m^{(1)} &= \text{SiLU}\left(W_m^{(1)} \mathbf{f}(z,A) + b_m^{(1)}\right), 
\quad W_m^{(1)} \in \mathbb{R}^{128 \times 33}, \\
\mathbf{h}_m^{(2)} &= \text{SiLU}\left(W_m^{(2)} \mathbf{h}_m^{(1)} + b_m^{(2)}\right), 
\quad W_m^{(2)} \in \mathbb{R}^{64 \times 128}.
\end{align}

The two pathways are subsequently concatenated
\begin{equation}
\mathbf{h}^{(3)} = \left[\mathbf{h}_c^{(2)}, \mathbf{h}_m^{(2)}\right] \in \mathbb{R}^{128}.
\end{equation}

This fused representation is processed through a residual block \cite{he2016resnet} to enhance
representational capacity while stabilizing optimization.
\begin{align}
\mathbf{h}^{(4)} &= \text{SiLU}\left(W^{(4)} \mathbf{h}^{(3)} + b^{(4)}\right), 
\quad W^{(4)} \in \mathbb{R}^{128 \times 128}, \\
\mathbf{h}^{(5)} &= \text{SiLU}\left(W^{(5)} \mathbf{h}^{(4)} + b^{(5)}\right), 
\quad W^{(5)} \in \mathbb{R}^{128 \times 128}, \\
\end{align}
The residual connection is then applied as
\begin{align}
\mathbf{h}^{(6)} &= \mathbf{h}^{(4)} + \mathbf{h}^{(5)}
\end{align}

The representation is then projected to a lower dimension
\begin{equation}
\mathbf{h}^{(7)} = \text{SiLU}\left(W^{(7)} \mathbf{h}^{(6)} + b^{(7)}\right), 
\quad W^{(7)} \in \mathbb{R}^{64 \times 128}.
\end{equation}

Finally, a sigmoid layer produces the probability that point $z$ lies within the $\varepsilon$-pseudospectrum:
\begin{equation}
p(z, A; \theta) = \sigma\left(w^T \mathbf{h}^{(7)} + b\right),
\end{equation}
where $\sigma(x) = \frac{1}{(1 + e^{-x})}$ is the sigmoid function, $w \in \mathbb{R}^{64}$, and $\theta$ denotes all trainable parameters.

We use the SiLU (Swish) activation function \cite{ramachandran2017searching} 
$\text{SiLU}(x) = x \cdot \sigma(x)$ throughout the network, which provides 
smooth, non-monotonic gradients that facilitate optimization. The network contains approximately 45,000 trainable parameters, making it lightweight enough for rapid inference while remaining sufficiently expressive to capture complex pseudospectral geometries. Having described the network architecture, we now detail the procedure used to construct the training data.

\subsection{Training Data Generation}

The training data are generated using random banded non-normal matrices, which provide a computationally efficient yet sufficiently rich class of matrices exhibiting nontrivial pseudospectral behavior. For each matrix $A \in \mathbb{R}^{64 \times 64}$, the bandwidth $\beta \in \{1,2,3,4\}$ is sampled uniformly, and the matrix entries are defined as
\begin{equation}\label{matrix}
A_{ij} = 
\begin{cases}
\text{uniform}\{-1, 0, 1\}, & \text{if } |i-j| \leq \beta, \\
0, & \text{otherwise}.
\end{cases}
\end{equation}
We additionally enforce $A \neq A^T$ to ensure non-normality and restrict the condition number to $\kappa(A) < 10^8$ to avoid extreme numerical instability.

For each matrix, ground-truth pseudospectra are computed on a $100 \times 100$ grid over the domain $[-4,4] \times [-4,4]$ in the complex plane with threshold $\varepsilon = 0.01$, which provides a suitable resolution for capturing relevant spectral sensitivity.

A grid point $z$ is labeled as sensitive if $\sigma_{\min}(zI - A) \leq \varepsilon$.

To address the inherent class imbalance between sensitive and non-sensitive grid points, we adopt the following sampling strategy.

\subsubsection{Balanced Sampling Strategy}

Direct grid sampling produces severe class imbalance, with typically fewer than $5\%$ of grid points classified as sensitive. To address this issue, we adopt a targeted sampling strategy that partitions the data into two classes \cite{he2009learning,buda2018systematic}. For each matrix, all sensitive grid points are included in the positive class. For the negative class, we randomly sample $\max(10\,n_{\text{pos}},\,200)$ non-sensitive points, where $n_{\text{pos}}$ denotes the number of sensitive points for that matrix. 

This strategy ensures sufficient representation of the positive class while maintaining a manageable dataset size. Using $500$ training matrices, this procedure yields a total of $438{,}084$ labeled samples. The resulting dataset is then used to optimize the network parameters as described below.

\subsection{Loss Function and Optimization}

We train the neural network using binary cross-entropy \cite{goodfellow2016deep}
\begin{equation}
\mathcal{L}(\theta) = -\frac{1}{N}\sum_{i=1}^N 
\left[ y_i \log p_i + (1-y_i) \log(1-p_i) \right],
\end{equation}
where $y_i \in \{0,1\}$ denotes the ground-truth label and
$p_i = p(z_i, A_i; \theta)$ is the predicted probability. Optimization is performed using the Adam optimizer \cite{kingma2014adam} with a learning rate of $10^{-3}$ and batch size $512.$ The network is trained for up to $25$ epochs with early stopping based on validation loss, using a patience of 5 epochs. 

During inference, evaluating the classifier efficiently over the complex plane becomes critical.

\subsection{Hierarchical Prediction Strategy}

Direct evaluation of the network on a $100 \times 100$ grid requires 10{,}000 forward passes. We mitigate this cost using a hierarchical coarse-to-fine prediction strategy \cite{berger1989local}. We first evaluate the network on a coarse $25 \times 25$ grid, corresponding to every fourth grid point. Let $P_{\text{coarse}}$ denote the resulting predicted probability map. Candidate regions are identified using an 80th-percentile threshold
\begin{equation}
\tau_{\text{coarse}} = \text{quantile}_{0.8}(P_{\text{coarse}}).
\end{equation}

For each coarse grid cell $[i:i+4,\, j:j+4]$ satisfying
$P_{\text{coarse}}[i,j] \ge \tau_{\text{coarse}}$, the network is evaluated on
all $4 \times 4 = 16$ fine-grid points within that cell. Cells below the
threshold are assigned zero probability. This hierarchical procedure reduces the number of network evaluations from 10{,}000 to approximately
$625 + 0.2\,(10{,}000 - 625) \approx 2{,}625$ on average, with minimal impact on prediction accuracy. To convert predicted probabilities into a binary sensitive region, an appropriate decision threshold is required.

\subsection{Threshold Calibration}

The raw network probabilities must be converted into binary predictions. To this end, we calibrate the classification threshold $\tau$ on a held-out validation set of 30 matrices, with the objective of maximizing recall while controlling false positives. For each candidate threshold $\tau \in \{0.05, 0.06, \ldots, 0.94\}$, the network predictions $\hat{S} = \{z : p(z, A) \geq \tau\}$ are first post-processed using morphological dilation with a $5 \times 5$ structuring element \cite{serra1982image, haralick1992computer} to provide a safety margin. Recall is then computed against the ground-truth sensitive regions, which are themselves dilated using a $3 \times 3$ structuring element to account for minor spatial misalignments. The final threshold is selected according to
\begin{equation}
\tau^* = \min \left\{ \tau :
\text{median}_{\text{val}}(\text{recall}) \geq 0.90 \;\text{and}\;
\text{P}_{10}(\text{recall}) \geq 0.75 \right\},
\end{equation}
where $\text{P}_{10}$ denotes the 10th percentile across the validation set. This criterion enforces both strong typical performance and robustness in worst-case scenarios. On the validation set, we obtain $\tau^* = 0.05$, yielding a median recall of $1.000$ and a 10th-percentile recall of $1.000$.

In practice, the proposed method integrates neural network guidance with
selective singular value computations in a coarse-to-fine framework. The network
is first used to localize a small subset of grid points likely to belong to the
$\varepsilon$-pseudospectrum, after which exact singular value decompositions are
performed only within this predicted region. The final pseudospectrum is then
constructed by thresholding $\sigma_{\min}(zI-A)$ on the restricted grid, while
all remaining points are classified as non-sensitive.

From a computational standpoint, the dominant cost arises from the singular
value decompositions. Network inference scales linearly with the number of
evaluated grid points and is carried out on a coarse grid followed by refinement
on only a limited fraction of fine-grid locations. This results in
$O(N_{\text{coarse}} + f\,N_{\text{fine}})$ network evaluations, where
$f \approx 0.26$ denotes the fraction of fine-grid points selected for
refinement. Exact singular value computations are subsequently restricted to the
predicted region, leading to a cost proportional to
$O(|x_{\text{in}}|\,|y_{\text{in}}|\,n^3)$, with
$|x_{\text{in}}|\,|y_{\text{in}}| \approx 0.16\,N_xN_y$ in our experiments.
As a result, evaluating only a small portion of the grid yields substantial
computational savings compared to full pseudospectrum computation.

The above analysis clarifies how restricting singular value decompositions to a carefully selected subset of the complex plane can substantially reduce computational cost without compromising accuracy. We now formalize this idea by presenting a hybrid computational framework that integrates neural network prediction with exact numerical evaluation. The goal is to translate the learned localization capability into a principled algorithmic strategy with transparent computational complexity.

\section{Hybrid Computational Approach}
\label{sec:hybrid_approach}

Having established the mathematical foundations in Sec.~\ref{sec:preliminaries} and detailed the neural network architecture in Sec.~\ref{sec:neural-approach}, we now describe the complete hybrid computational strategy. Given a matrix $A \in \mathbb{C}^{n \times n}$ and a computational grid $\mathcal{G}$ defined over a rectangular domain in the complex plane, our approach computes an approximation $\widetilde{\mathcal{S}}_\varepsilon$ to the sensitive zone $\mathcal{S}_\varepsilon$ through a two-stage process.

First, a trained classifier $f_\theta : \mathbb{C} \times \mathbb{R}^{n \times n} \to [0,1]$ identifies a candidate region $\mathcal{R} \subseteq \mathcal{G}$ that likely contains $\mathcal{S}_\varepsilon$, with $|\mathcal{R}| \ll N$ to achieve computational savings. The classifier is designed to satisfy 
\[
\mathbb{P}[\mathcal{S}_\varepsilon \subseteq \mathcal{R}] \geq 0.90,
\]
through threshold calibration on validation data, where $\mathbb{P}$ denotes probability with respect to the distribution of matrices and grid points used during validation.

Second, exact singular value decompositions are computed only for points $z \in \mathcal{R}$, yielding $\sigma_{\min}(z) = \sigma_{\min}(zI - A)$ to machine precision. The approximate sensitive zone is then defined as 
\[
\widetilde{\mathcal{S}}_\varepsilon = \{z \in \mathcal{R} : \sigma_{\min}(z) \leq \varepsilon\}.
\]

This two-stage procedure yields the hybrid computational cost
\begin{equation}
\label{eq:hybrid_cost}
\mathcal{C}_{\text{hybrid}} = \mathcal{C}_{\text{NN}}(N) + |\mathcal{R}| \cdot \mathcal{C}_{\text{SVD}}(n),
\end{equation}
where $\mathcal{C}_{\text{NN}}(N) = O(N)$ represents the cost of neural network inference over the grid and $\mathcal{C}_{\text{SVD}}(n) = O(n^3)$ is the cost of a single SVD.

Since neural network evaluation is substantially cheaper than SVD—typically $\mathcal{C}_{\text{NN}}(N) \approx 10^{-4} \cdot N \cdot \mathcal{C}_{\text{SVD}}(n)$ for moderate $n$—and the predicted region satisfies $|\mathcal{R}| \ll N$ by construction, equation~\eqref{eq:hybrid_cost} represents a significant reduction compared to the full cost $\mathcal{C}_{\text{full}} = N \cdot \mathcal{C}_{\text{SVD}}(n)$ from equation~\eqref{eq:full_cost}.

The theoretical speedup is approximately
\begin{equation}\label{eq:speedup_theoretical}
\text{Speedup}_{\text{theory}}
\approx \frac{\mathcal{C}_{\text{full}}}{\mathcal{C}_{\text{hybrid}}}
\approx \frac{N}{|\mathcal{R}|}
= \rho^{-1},
\end{equation}
where $\rho = |\mathcal{R}|/N$ denotes the fraction of the grid requiring SVD computation. In our experiments with banded non-normal matrices, we observe $\rho \approx 0.16$, suggesting potential speedups of approximately $6\times$.

To minimize the cost of neural network evaluations themselves, we employ a hierarchical coarse-to-fine prediction strategy. Rather than evaluating $f_\theta$ at all $N$ grid points, we first construct a coarse grid $\mathcal{G}_{\text{coarse}}$ by subsampling $\mathcal{G}$ with stride $s = 4$ in both spatial dimensions, yielding $N_{\text{coarse}} = N/s^2$ points. The classifier is evaluated on $\mathcal{G}_{\text{coarse}}$, and high-probability regions are identified via adaptive thresholding: $\mathcal{R}_{\text{coarse}} = \{z \in \mathcal{G}_{\text{coarse}} : p(z) \geq \tau_{\text{coarse}}\}$, where $\tau_{\text{coarse}}$ is defined as the 80th percentile of coarse-grid predictions $\{p(z) : z \in \mathcal{G}_{\text{coarse}}\}$. This percentile-based threshold adapts to the specific pseudospectral structure of each matrix. 

All fine-grid points lying within cells flagged by $\mathcal{R}_{\text{coarse}}$ are then collected into a refinement region $\mathcal{R}_{\text{fine}}$, and the classifier is evaluated on this subset. A calibrated decision threshold $\tau$—determined independently via the validation procedure described in Section~\ref{sec:neural-approach}—is applied to obtain the final predicted region $\mathcal{R} = \{z \in \mathcal{R}_{\text{fine}} : p(z) \geq \tau\}$. To ensure robustness near pseudospectral boundaries and account for discretization effects, $\mathcal{R}$ is expanded using a $5 \times 5$ morphological dilation. This hierarchical scheme reduces the total number of network evaluations from $N$ to approximately $N/16 + |\mathcal{R}_{\text{fine}}| \approx 2{,}625$ for a $100 \times 100$ grid, where typically $|\mathcal{R}_{\text{fine}}| \approx 0.2N$. 

\subsection{Complexity Analysis}

The total computational cost decomposes as
\begin{equation}
\mathcal{C}_{\text{hybrid}} = \underbrace{O(n^3)}_{\text{eigenvalues}} + \underbrace{O(N/s^2)}_{\text{coarse NN}} + \underbrace{O(|\mathcal{R}_{\text{fine}}|)}_{\text{fine NN}} + \underbrace{|\mathcal{R}| \cdot O(n^3)}_{\text{restricted SVD}}.
\end{equation}

For typical problem sizes ($n = 64$, $N = 10^4$, $|\mathcal{R}| \approx 0.16N$), the dominant term is the restricted SVD stage. The actual speedup realized in practice is
\begin{equation}
\text{Speedup}_{\text{actual}} = \frac{t_{\text{full}}}{t_{\text{NN}} + t_{\text{restricted SVD}}},
\end{equation}
where timing measurements account for implementation overhead. As shown in Section~\ref{sec:experiments}, we observe speedups of $2.45\times$ on average, with the discrepancy from theoretical speedup~\eqref{eq:speedup_theoretical} arising from neural network inference time and conservative region padding.

\subsection{Error Characteristics}

The proposed hybrid strategy admits only a limited and well-controlled form of approximation error. If $\mathcal{S}_\varepsilon \not\subseteq \mathcal{R}$, some truly sensitive points are omitted. This is controlled via threshold calibration targeting $\geq 90\%$ recall on validation data. Critically, the method computes \emph{exact} $\sigma_{\min}(z)$ for all $z \in \mathcal{R}$ via SVD. Thus, within the predicted region, all values are correct to machine precision. Errors manifest only as omissions, not incorrect predictions.

This error structure is favorable for applications where conservative estimates are acceptable, provided that computed values are reliable where available. We now empirically assess the effectiveness of this approach in terms of both accuracy and computational efficiency through a series of numerical experiments.

\section{Numerical Experiments}\label{sec:experiments}

We evaluate the proposed hybrid pseudospectra computation method on randomly generated banded non-normal matrices and compare its performance against both exhaustive grid evaluation and a random sampling baseline.

\subsection{Experimental Setup}

Test matrices are generated from the same class of random banded non-normal
matrices defined in \eqref{matrix}, with matrix size $A \in \mathbb{R}^{64 \times 64}$ and bandwidth $\beta \in \{1,2,3,4\}$ sampled uniformly. This ensures that the test problems share the same structural characteristics as the training data while remaining disjoint from the training and validation sets. The neural network is trained using the procedure described in Sec.~\ref{sec:neural-approach} and is kept fixed for all experiments reported in this section. Threshold calibration is performed on a separate validation set of $30$ matrices, yielding a calibrated decision threshold $\tau^* = 0.05$, which is used consistently throughout all numerical experiments.

Using this trained model and calibrated threshold, we evaluate the performance on 50 held-out test matrices generated from the
same distribution as the training data but disjoint from both the training and validation sets. For each test matrix, the ground-truth
$\varepsilon$-pseudospectrum is computed using full singular value decomposition on a $100 \times 100$ grid over $[-4,4] \times [-4,4]$ in the complex plane with $\varepsilon = 0.01$. The proposed hybrid method combines neural network predictions with restricted singular value decompositions in a coarse-to-fine manner, as described in the preceding section. For comparison, we also consider a random sampling baseline that evaluates an identical fraction of grid points selected uniformly at random. All timing measurements are performed using a single CPU thread in order to exclude any effects of parallelization or hardware acceleration. The accuracy and efficiency of these evaluations are quantified using the
performance metrics described below.

\subsection{Performance Metrics}

We assess the proposed method using both classification accuracy measures and computational efficiency indicators. To account for boundary uncertainty in the pseudospectrum, the ground-truth sensitive region is dilated using a $3 \times 3$ structuring element prior to metric computation.Let $S_{\text{true}}^{\text{dilated}}$ denote the dilated ground-truth region and $S_{\text{true}}$ the original sensitive region.

Classification performance is quantified by accuracy, precision, recall, and coverage. Accuracy measures the fraction of correctly classified grid points. Precision is defined as 
\[
\frac{|S_{\text{true}} \cap S_{\text{pred}}|}{|S_{\text{pred}}|},
\]
reflecting the purity of the predicted sensitive region. 

Recall is defined as 
\[
\frac{|S_{\text{true}}^{\text{dilated}} \cap S_{\text{pred}}|}{|S_{\text{true}}^{\text{dilated}}|},
\]
which measures the fraction of true sensitive points successfully identified. Recall is computed with respect to the dilated ground-truth region, while coverage is evaluated on the original sensitive region, providing a complementary measure of prediction completeness.

Coverage is defined as the fraction of true sensitive points captured by the predicted region, i.e.,
\[
\frac{|S_{\text{true}} \cap S_{\text{pred}}|}{|S_{\text{true}}|}.
\]

Computational efficiency is evaluated by measuring the fraction of grid points on which singular value decompositions are performed, referred to as the grid fraction. We further report the actual speedup, defined as
$\frac{t_{\text{full}} }{(t_{\text{NN}} + t_{\text{restricted}})}$, which accounts for both neural network inference and restricted singular value computations, as well as the best-case speedup, $\frac{t_{\text{full}} }{ t_{\text{restricted}}}$, which excludes neural network overhead and represents an upper bound on achievable acceleration. 

Having defined the evaluation criteria, we now summarize the empirical
performance of the proposed hybrid method across a collection of unseen test matrices. The following results quantify the trade-off between accuracy and computational efficiency achieved by restricting expensive singular value computations to network-predicted regions of interest.

\subsection{Main Results}

Table~\ref{tab:main_results} summarizes performance across 50 test matrices. The method achieves high recall ($99.5\% \pm 1.4\%$) and coverage ($99.8\% \pm 0.7\%$) while evaluating only $15.9\% \pm 3.2\%$ of grid points, resulting in $2.45\times$ average speedup over exhaustive evaluation.

\begin{table}[t]
\centering
\caption{Performance statistics over 50 test matrices (mean $\pm$ std).}
\label{tab:main_results}
\begin{tabular}{lrrrr}
\toprule
\textbf{Metric} & \textbf{Mean} & \textbf{Std} & \textbf{Median} & \textbf{Range} \\
\midrule
\multicolumn{5}{l}{\textit{Classification Performance}} \\
Accuracy & 0.869 & 0.024 & 0.870 & [0.815, 0.927] \\
Precision & 0.175 & 0.059 & 0.173 & [0.058, 0.284] \\
Recall & 0.995 & 0.014 & 1.000 & [0.922, 1.000] \\
Coverage & 0.998 & 0.007 & 1.000 & [0.964, 1.000] \\
\midrule
\multicolumn{5}{l}{\textit{Computational Efficiency}} \\
Grid fraction & 0.159 & 0.032 & 0.166 & [0.090, 0.219] \\
Speedup (actual) & 2.45 & 0.36 & 2.50 & [1.50, 3.09] \\
Speedup (best-case) & 3.52 & 0.91 & 3.34 & [1.77, 5.84] \\
\midrule
\multicolumn{5}{l}{\textit{Timing (seconds)}} \\
Full SVD & 3.73 & 0.99 & 4.14 & [1.73, 4.85] \\
NN prediction & 0.38 & 0.07 & 0.36 & [0.31, 0.66] \\
Restricted SVD & 1.19 & 0.52 & 1.32 & [0.32, 2.39] \\
Hybrid total & 1.57 & 0.53 & 1.65 & [0.64, 2.82] \\
\bottomrule
\end{tabular}
\end{table}

The low precision (17.5\%) reflects conservative prediction: the network identifies approximately $16\%$ of the grid as potentially sensitive, while ground truth typically contains $<1\%$ sensitive points (mean 83 points out of 10,000). However, high recall and coverage confirm that nearly all true sensitive points are captured within predicted regions.

\begin{figure}[t]
\centering
\includegraphics[width=\textwidth]{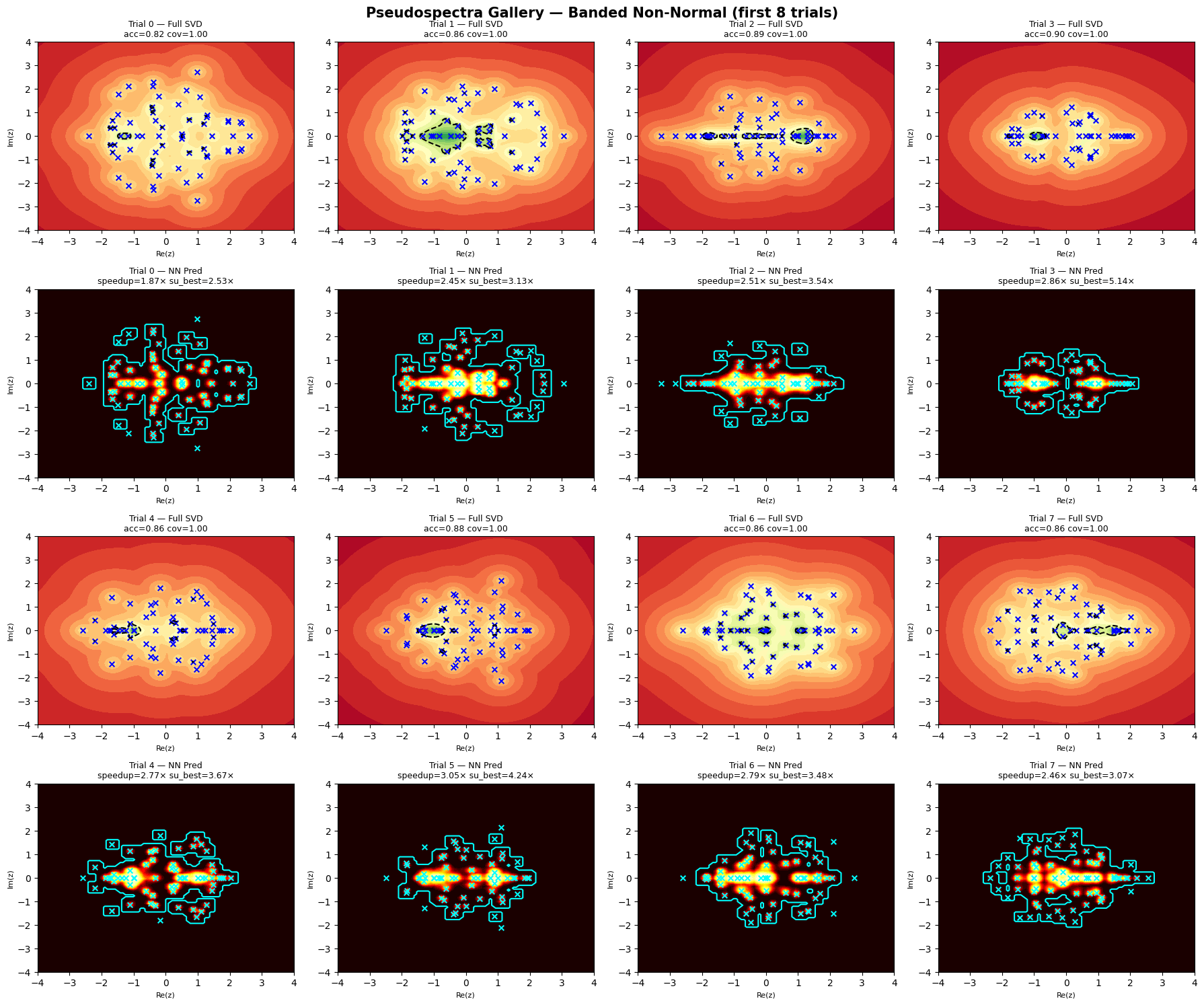}
\caption{Pseudospectra comparison on eight representative test matrices. For each matrix, top panel shows ground truth $\log(\sigma_{\min}(zI-A))$ with $\varepsilon$-pseudospectra boundary (black contour) and eigenvalues (blue crosses). Bottom panel shows neural network probability map with predicted sensitive region (cyan contour) and eigenvalues (cyan crosses). The method accurately captures diverse pseudospectra geometries across bandwidth values $\beta \in \{1,2,3,4\}$.}
\label{fig:gallery}
\end{figure}

\begin{figure}[t]
\centering
\includegraphics[width=\textwidth]{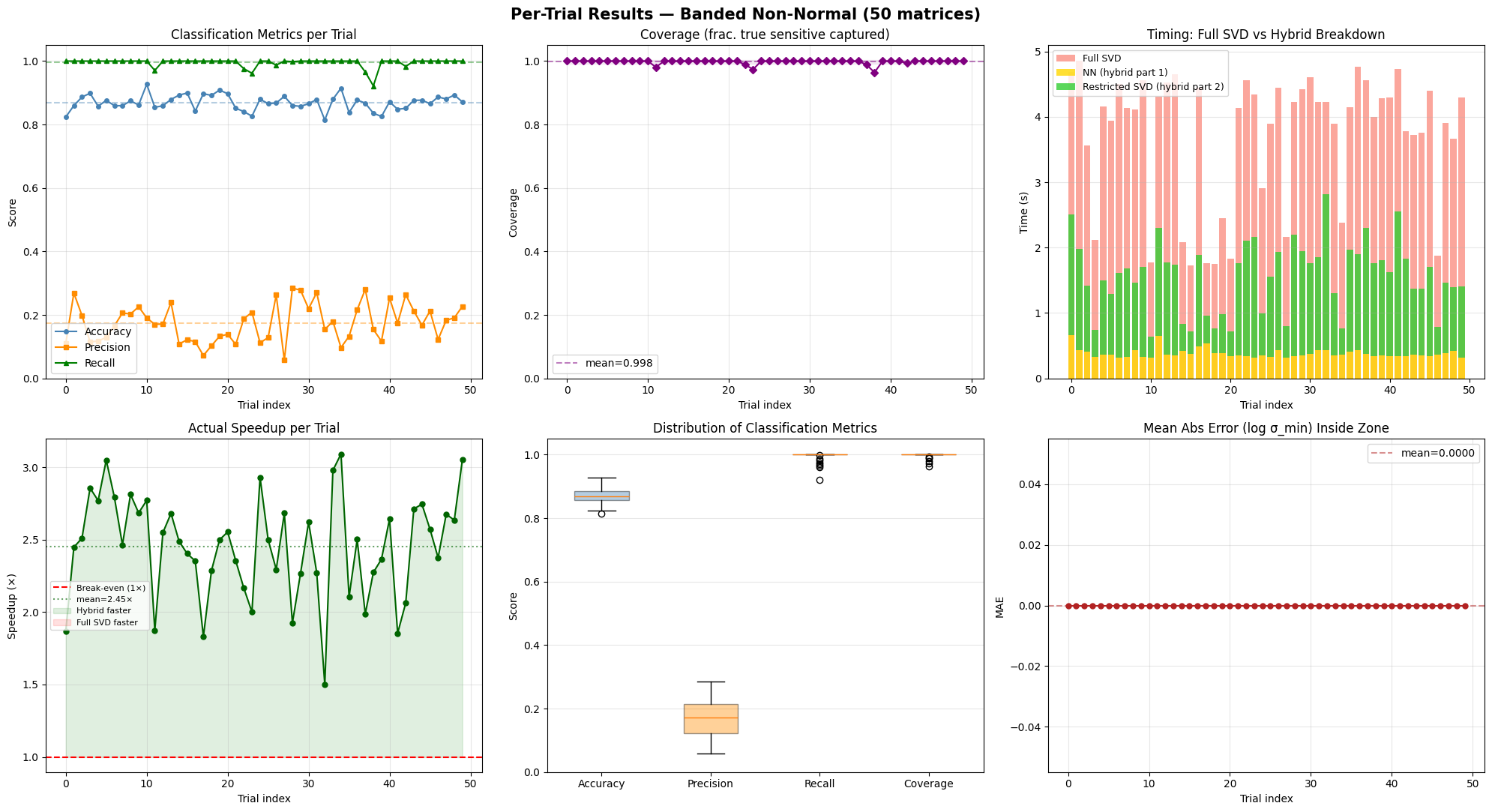}
\caption{Performance metrics across 50 test matrices. (a) Classification metrics per trial: accuracy (blue circles), precision (orange squares), and recall (green triangles) with mean values shown as dashed lines. (b) Coverage per trial (purple diamonds) showing fraction of true sensitive points captured. (c) Timing breakdown: full SVD (salmon), neural network overhead (gold), and restricted SVD (green). (d) Speedup per trial showing all matrices achieve $>1.5\times$ acceleration. (e) Distribution of classification metrics. (f) Mean absolute error in $\log(\sigma_{\min})$ within predicted zones.}
\label{fig:performance}
\end{figure}

Figure~\ref{fig:gallery} illustrates representative examples showing ground truth pseudospectra alongside neural network predictions. The network successfully identifies sensitive regions of varying shapes and extents across different matrix structures. Figure~\ref{fig:performance} shows classification metrics and speedup remain consistently high across all 50 test matrices, with the hybrid method achieving faster computation than full SVD in every trial.

To further examine how structural properties of the matrix influence performance, we next analyze the results as a function of matrix bandwidth, which directly controls sparsity and the complexity of the resulting pseudospectra.

\subsection{Bandwidth Stratification}

Table~\ref{tab:bandwidth} stratifies results by matrix bandwidth. Performance degrades slightly for larger bandwidth: matrices with $\beta = 1$ achieve $2.56\times$ speedup with perfect recall, while $\beta = 2$ achieves $2.70\times$ speedup. For larger bandwidth values, performance gradually degrades: $\beta = 4$ yields $2.06\times$ speedup with $98.4\%$ recall. This occurs because larger bandwidth produces more complex pseudospectra with broader sensitive regions, reducing the benefit of domain restriction.

\begin{table}[t]
\centering
\caption{Performance stratified by bandwidth (mean $\pm$ std).}
\label{tab:bandwidth}
\begin{tabular}{lcccc}
\toprule
\textbf{Bandwidth} & \textbf{Count} & \textbf{Speedup} & \textbf{Recall (\%)} & \textbf{Coverage (\%)} \\
\midrule
$\beta = 1$ & 12 & $2.56 \pm 0.32$ & $100.0 \pm 0.0$ & $100.0 \pm 0.0$ \\
$\beta = 2$ & 11 & $2.70 \pm 0.20$ & $100.0 \pm 0.0$ & $100.0 \pm 0.0$ \\
$\beta = 3$ & 15 & $2.50 \pm 0.28$ & $99.7 \pm 0.8$ & $99.9 \pm 0.5$ \\
$\beta = 4$ & 12 & $2.06 \pm 0.25$ & $98.4 \pm 2.3$ & $99.3 \pm 1.1$ \\
\bottomrule
\end{tabular}
\end{table}

Despite this variation, the method maintains $>98\%$ recall across all bandwidth values, demonstrating robustness to structural diversity within the banded matrix class. To further contextualize these results, we compare the proposed method with a random sampling baseline that evaluates the same fraction of grid points but without exploiting geometric information.

\subsection{Comparison with Random Sampling}

To isolate the value of geometric guidance, we compare against a random sampling baseline that evaluates the same fraction of grid points ($15.9\%$) but selects them uniformly at random rather than using neural network predictions.

Table~\ref{tab:comparison} shows that random sampling achieves only $45.9\% \pm 14.2\%$ recall compared to our method's $99.5\% \pm 1.4\%$ recall. The high standard deviation ($14.2\%$) in random sampling indicates unreliable performance: some matrices yield acceptable coverage by chance, while others miss most sensitive points entirely.

\begin{table}[t]
\centering
\caption{Method comparison (mean $\pm$ std over 50 matrices).}
\label{tab:comparison}
\begin{tabular}{lcccc}
\toprule
\textbf{Method} & \textbf{Coverage (\%)} & \textbf{Recall (\%)} & \textbf{Precision (\%)} & \textbf{Speedup} \\
\midrule
Full SVD & 100.0 & 100.0 & -- & $1.00\times$ \\
Random sampling & $15.9 \pm 3.2$ & $45.9 \pm 14.2$ & $66.7 \pm 10.9$ & $6.28\times$ (theor.) \\
Hybrid (ours) & $15.9 \pm 3.2$ & $99.5 \pm 1.4$ & $17.5 \pm 5.9$ & $2.45\times$ (actual) \\
\bottomrule
\end{tabular}
\end{table}

\begin{figure}[t]
\centering
\includegraphics[width=\textwidth]{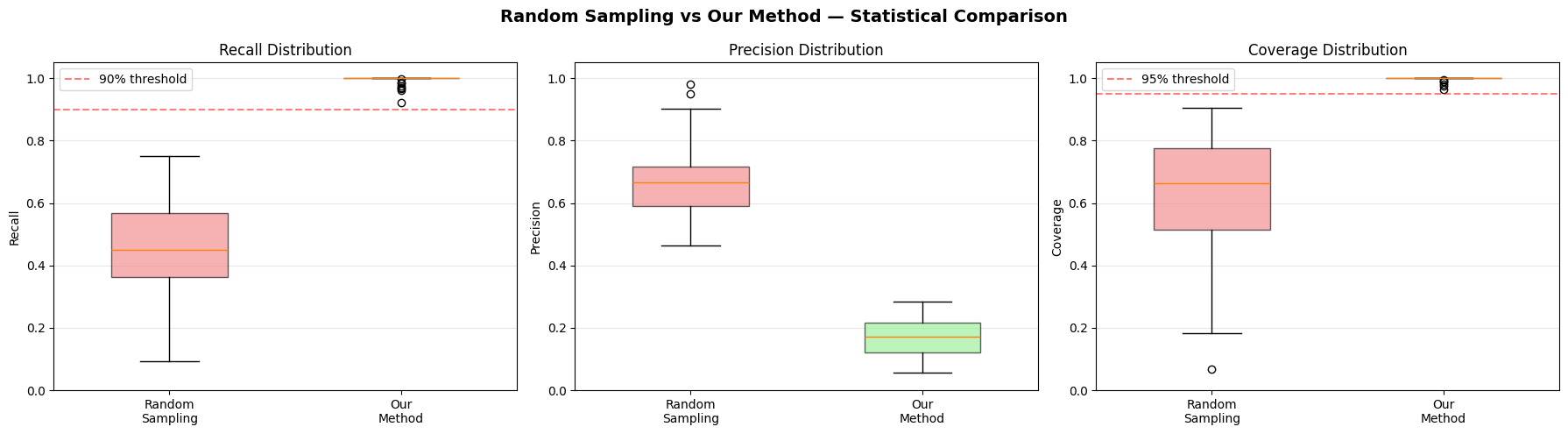}
\caption{Comparison with random sampling baseline over 50 test matrices. Box plots show (a) recall, (b) precision, and (c) coverage distributions. Random sampling (red) exhibits high variance and poor median recall (45.9\%), while our method (green) achieves consistent high recall (99.5\%) with low variance. Red dashed lines indicate target thresholds.}
\label{fig:baseline}
\end{figure}

Figure~\ref{fig:baseline} compares the distributions of recall, precision, and coverage between random sampling and our method. The dramatic improvement in recall consistency ($1.4\%$ vs $14.2\%$ standard deviation) demonstrates that learned geometric structure drives reliable performance.

The theoretical speedup for random sampling ($6.28\times$) assumes zero overhead, while our actual speedup ($2.45\times$) includes neural network evaluation time ($0.38$s average). However, random sampling's poor recall renders the theoretical speedup meaningless: a method that misses half the sensitive points cannot reliably compute pseudospectra.

Our method achieves $39.1\%$ efficiency relative to the theoretical maximum
($2.45 / 6.28 = 0.391$), trading some speed for substantially improved
reliability. The $+53.6$ percentage point gain in recall over random sampling demonstrates that accuracy is driven by learned geometric structure rather than mere reduction in computational budget. We now examine the computational efficiency of the proposed approach by analyzing its runtime characteristics and resulting speedup behavior.

\subsection{Timing Analysis}

\begin{figure}[t]
\centering
\includegraphics[width=0.9\textwidth]{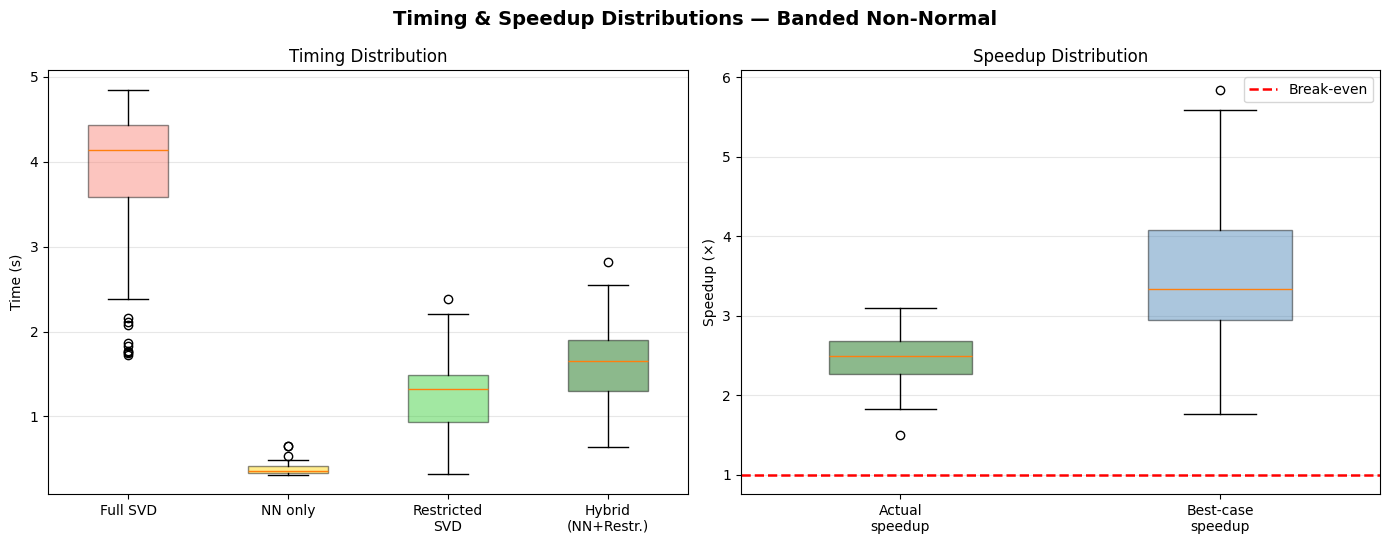}
\caption{Timing and speedup distributions over 50 test matrices. Left panel: absolute timing for full SVD (red), neural network only (yellow), restricted SVD (green), and total hybrid method (dark green). Right panel: actual speedup (green) and best-case speedup excluding network overhead (blue), with break-even line at $1\times$ (red dashed).}
\label{fig:timing}
\end{figure}

Figure~\ref{fig:timing} breaks down timing components and speedup distributions. The neural network overhead ($0.38$s average) is modest relative to the time saved by restricted evaluation: full SVD requires $3.73$s average, while our hybrid approach requires $1.57$s average.

The difference between the observed speedup ($2.45\times$) and the best-case speedup ($3.52\times$) is attributable to the overhead introduced by neural network inference. As matrix dimension increases and the cost of singular value decomposition scales as $O(n^3)$, this overhead becomes progressively less significant, indicating the potential for greater acceleration at larger problem sizes. Network evaluation itself scales sublinearly with grid size due to the hierarchical prediction strategy: the coarse stage requires only $625$ evaluations, and refinement is performed on approximately $26\%$ of the fine grid, resulting in an average of $2{,}625$ network evaluations compared to $10{,}000$ for full-grid evaluation. This hierarchical design substantially reduces inference overhead while preserving prediction accuracy.

\section{Discussion}\label{sec:discussion}

The effectiveness of the proposed hybrid framework arises from the combination of three key components: geometry-aware feature design, conservative threshold calibration, and hierarchical prediction. Together, these elements enable reliable identification of pseudospectrally sensitive regions while significantly reducing the number of singular value decompositions.

The feature representation integrates global matrix descriptors with local eigenvalue distance information, allowing the model to capture both overall pseudospectral behavior and spatial variation. In addition, features related to eigenvector conditioning and resolvent estimates provide indicators of spectral sensitivity beyond eigenvalue location alone.

The method is intentionally conservative: threshold calibration and morphological dilation prioritize high recall, ensuring that sensitive regions are not missed. This design reflects the asymmetric cost of errors, where false negatives are significantly more harmful than false positives.

The hierarchical coarse-to-fine strategy further improves efficiency by restricting both neural network evaluation and SVD computations to candidate regions, reducing overhead without sacrificing accuracy.

A limitation of the current study is its focus on banded non-normal matrices. While the framework is general, its performance on other matrix classes remains to be investigated. Moreover, training requires labeled pseudospectra, which introduces an upfront computational cost.

Overall, the results demonstrate that learned geometric structure can effectively guide numerical computation, enabling substantial acceleration while maintaining reliability.

\section{Conclusion}\label{sec:conclusion}

We presented a neural-guided hybrid framework for accelerating pseudospectra computation of non-normal matrices. By predicting spectrally sensitive regions using matrix features, the method restricts expensive singular value decompositions to a small subset of the complex plane.

Experiments on banded non-normal matrices show that the approach evaluates only about $16\%$ of grid points while achieving approximately $2.45\times$ speedup and maintaining high reliability, with $99.5\%$ recall and $99.8\%$ coverage.

The results demonstrate that data-driven domain restriction can effectively complement classical numerical methods, providing a practical strategy for efficient pseudospectra computation. Future work includes extending the approach to broader matrix classes and evaluating performance at larger problem scales.

\bibliographystyle{unsrt}
\bibliography{references}

\end{document}